\documentclass{ifacconf}

\usepackage{graphicx}      
\usepackage{subcaption}
\usepackage{amsmath}

\usepackage{dsfont}
\usepackage{xcolor}
\usepackage{upgreek}
\usepackage{enumerate}
\usepackage{amsfonts}
\usepackage{amssymb}
\usepackage{textcomp}
\usepackage{natbib}        

\usepackage{tikz}
\usepackage{mathabx}
\usepackage{tkz-berge}
\usetikzlibrary{calc}
\usetikzlibrary{shapes,fit}
\usepackage{verbatim}
\usepackage{setspace}
\usetikzlibrary{shapes,arrows}

\newtheorem{theorem}{Theorem}[section]

\newtheorem{lemma}[theorem]{Lemma}
\newtheorem{assumption}[theorem]{Assumption}

\newtheorem{corollary}[theorem]{Corollary}

\newtheorem{remark}[theorem]{Remark}

\DeclareMathOperator{\spa}{span}

\DeclareMathOperator{\sign}{sign}

\DeclareMathOperator{\diag}{diag}

\newcommand{\co} {\ensuremath{\mathrm{co}}} 

\newcommand{\calC}{\ensuremath{\mathcal{C}}}
\newcommand{\calD}{\ensuremath{\mathcal{D}}}
\newcommand{\calE}{\ensuremath{\mathcal{E}}}
\newcommand{\calF}{\ensuremath{\mathcal{F}}}
\newcommand{\calG}{\ensuremath{\mathcal{G}}}

\newcommand{\calI}{\ensuremath{\mathcal{I}}}

\newcommand{\calK}{\ensuremath{\mathcal{K}}}
\newcommand{\calL}{\ensuremath{\mathcal{L}}}

\newcommand{\calQ}{\ensuremath{\mathcal{Q}}}
\newcommand{\calR}{\ensuremath{\mathcal{R}}}

\newcommand{\calV}{\ensuremath{\mathcal{V}}}

\newcommand{\R}{\ensuremath{\mathbb R}}

\newcommand{\Z}{\ensuremath{\mathbb Z}}

\newcommand{\funcRdR}{\ensuremath{{f}}}
\newcommand{\funcRdRd}{\ensuremath{{X}}}

\begin{document}
\begin{frontmatter}

\title{ Nonlinear consensus protocols with applications to quantized systems \thanksref{footnoteinfo}}

\thanks[footnoteinfo]{This work was supported by the
Knut and Alice Wallenberg Foundation, the  Swedish Foundation for Strategic Research, and the Swedish Research Council.}

\author{Jieqiang Wei, Xinlei Yi, Henrik Sandberg and Karl Henrik Johansson}

\address{ACCESS Linnaeus Centre and Electrical Engineering, \\KTH Royal Institute of Technology, 100 44, Stockholm, Sweden \\(e-mail: jieqiang, xinleiy, hsan, kallej@kth.se).}

\begin{abstract}\label{s:Abstract}
This paper studies multi-agent systems with nonlinear consensus protocols, i.e., only nonlinear measurements of the states are available to agents. The solutions of these systems are understood in Filippov sense since the possible discontinuity of the nonlinear controllers. Under the condition that the nonlinear functions are monotonic increasing without any continuous constraints, asymptotic stability is derived for systems defines on both directed and undirected graphs. The results can be applied to quantized consensus which extend some existing results from undirected graphs to directed ones.
\end{abstract}

\end{frontmatter}

\section{Introduction}
Distributed consensus a fundamental problem in the study of multi-agent systems.
In addition to the well studied linear consensus problem (see e.g., \cite{Olfati2004}, \cite{Moreau2004}, \cite{Ren2005}), the nonlinear version has attracted much attention. Generally speaking, for continuous time models, nonlinear consensus studies can be divided into continuous and discontinuous systems. For the continuous case, we refer to \cite{Papachristodoulou2010}, \cite{Lin2007}, \cite{Andreasson2012} etc.
In this paper, instead we focus on the nonlinear consensus protocol with discontinuous dynamics. There are several existing works about this topic. Here we review some of the most related ones.
In \cite{Cortes2006}, the author studied the finite-time convergence of
\begin{equation}\label{e:intro_sign}
\dot{x} =\sign(-Lx),
\end{equation}
where $L$ is the Laplacian matrix of the graph and $\sign$ is the signum function.
It is proved that Filippov solutions will converge to average-max-min consensus in finite time. However, the result is not precise in the sense that it does not hold for all solutions. In \cite{Wei2015}, the authors considered the more general model
\begin{equation}\label{e:intro_old_general}
\dot{x}=f(-Lx),
\end{equation}
where $f$ is any sign-preserving function, i.e., each component of $f$ takes positive value for positive argument and vice versa. Sufficient conditions to guarantee asymptotic consensus of  all Filippov solutions are given in \cite{Wei2015}. In \cite{Kashyap2007}, the authors considered a  discretized version of \eqref{e:intro_old_general} with $f$ being a quantizer and $L$ a time-varying stochastic matrix.

Motivated by some practical scenarios, such as multi-robot coordination with coarse measurements, the model to be investigated in this paper is
\begin{equation}\label{e:intro_nonlinear_measurement}
\dot{x}=-Lf(x),
\end{equation}
where we assume $f$ to be any monotone function not necessarily cross the origin. The measurement of the state of each agent can obey different nonlinear criteria: quantized, biased etc. One closely related existing work is \cite{Liu2015consensus},  where the authors employ a stronger assumption, i.e., the nonlinear function $f$ is piecewise continuous, strictly monotone and sign preserving. In \cite{Liu2015consensus}, precise consensus can be achieved. However, their stronger assumption puts limits on the applicability of the results, for example, quantized measurement maps fail to be strictly monotone. A special case of the system \eqref{e:intro_nonlinear_measurement} is $f$ equal to the uniform quantizer. For such systems, \cite{Ceragioli2011} and \cite{Frasca2012} showed the asymptotic convergence of all the Krasovskii solutions to practical consensus. Furthermore, they assume undirected graphs. We extend these results to directed cases.
For the system \eqref{e:intro_nonlinear_measurement}, we address the stability using the notion of Filippov solution.
The reasons we choose Filippov solution are following.  First, for many nonlinear consensus protocols with discontinuous controllers, the classical and Carath\'{e}odory solutions do not exist. For example, in \cite{Ceragioli2011}, it is proven that both classical and Carath\'{e}odory solutions do not exist in general for system \eqref{e:intro_nonlinear_measurement} with $f$ being a uniform quantizer. So considering generalized solutions is necessary. 
Second, Filippov solution, comparing to Krasovskii solution, can eliminate the irregular behavior from the general nonlinear differential inclusion. Third,  for  quantized systems, Filippov  and Krasovskii solutions are equivalent.

The contributions of this paper are twofolds. First, we consider the general nonlinear consensus protocol \eqref{e:intro_nonlinear_measurement}, and present a stability analysis for all Filippov solutions under the weakest fixed topology, namely directed graphs containing spanning trees. Our result incorporates many existing works as special cases. 
Second, we consider the special case of quantized consensus protocols and present an extension to the results in \cite{Ceragioli2011}, \cite{Frasca2012} from undirected graphs to directed ones.

The structure of the paper is as follows. In Section \ref{s:Preliminaries}, we introduce some preliminaries. 
In Section \ref{s:measure}, we prove convergence for nonlinear consensus protocols where the measurements of the state are effected by nonlinearities.
In Section \ref{s:quantize} we apply the results in Section \ref{s:measure} to quantized consensus protocols.
Finally, the paper is wrapped up with the conclusion in Section \ref{s:conclusion}.

\section{Preliminaries}\label{s:Preliminaries}
In this section we first briefly review some notions from graph
theory, e.g, \cite{Bollobas98, biggs1993algebraic}, and then give some properties of Filippov solutions (\cite{cortes2008}).

Let $\mathcal{G}=(\mathcal{V},\mathcal{E},A)$ be a weighted digraph with
node set $\mathcal{V}=\{v_1,\ldots,v_n\}$,
edge set $\mathcal{E}\subseteq\mathcal{V}\times\mathcal{V}$ and
weighted adjacency matrix $A=[a_{ij}]$ with nonnegative adjacency elements $a_{ij}$.
An edge of $\mathcal{G}$ is denoted by $e_{ij}=(v_i,v_j)$ and we write $\calI=\{1,2,\ldots,n\}$.
The adjacency elements $a_{ij}$ are associated with the edges of the graph in
the following way: $a_{ij}>0$  if and only if  $e_{ji} \in \mathcal{E}$.
Moreover,  $a_{ii}=0$ for all $i \in\calI$.
For undirected graphs,  $A=A^T$.

The set of neighbors of node $v_i$ is
denoted by $N_i=\{v_j\in\mathcal{V}:(v_j,v_i)\in\mathcal{E}\}$.
For each node $v_i$, its in-degree
is defined as
\begin{align*}
\deg_{\rm in} (v_i) = \sum_{j=1}^n a_{ij}.  
\end{align*}
The degree
matrix of the digraph $\mathcal{G}$ is a diagonal matrix $\Delta$
where $\Delta_{ii}=\deg_{\rm  in}(v_i)$. The \emph{graph Laplacian} is
defined as
\begin{equation*}
L=\Delta-A.
\end{equation*}
This implies $L\mathds{1}_n =0_n$, where $\mathds{1}_n$ is the $n$-vector containing only ones and $0_n$ is the $n$-vector containing only zeros.

A directed path from node $v_i$ to node $v_j$ is a chain of edges from $\calE$
such that the first edge starts from $v_i$, the last edge ends at  $v_j$ and
every edge starts where the previous edge ends.
A graph is called \emph{strongly connected} if for every two nodes $v_i$ and
$v_j$ there is a directed path from $v_i$ to $v_j$. A directed graph is called \emph{weakly connected} if replacing all of its directed edges with undirected edges produces a connected (undirected) graph.
A subgraph $\calG' = (\calV',\calE',A')$ of $\calG$ is called a \emph{directed spanning tree} for $\calG$ if $\calG'$ is weakly connected, $\calV' =\calV $, $\calE' \subseteq \calE$, $|\calE'|=n-1$, and for every node $v_i\in
\calV'$ there is exactly one $v_j$ such that $e_{ji}\in \calE'$, except for one
node, which is called the root of the spanning tree.
Furthermore, we call a node $v\in \calV$ a \emph{root} of $\calG$ if there is a directed
spanning tree for $\calG$ with $v$ as a root.
In other words, if $v$ is a root of $\calG$, then there is a directed path from
$v$ to every other node in the graph.

A digraph, with $m$ edges, is completely specified by its \emph{incidence matrix} $B$, which is an $n\times m$ matrix, with
element $(i,j)$ equal to $-1$ if the $j^{\text{th}}$ edge is towards vertex $i$, and equal to $1$ if the $j^{\text{th}}$ edge is originating from vertex $i$, and $0$ otherwise.


\begin{lemma} [\cite{Lu2007}] \label{connected_laplacian}
The graph Laplacian matrix $L$ of a strongly connected digraph
$\calG$ satisfies that zero is an algebraically simple
eigenvalue of $L$ and there is a positive vector
$w^{\top}=[w_{1},\cdots,w_{n}]$ such that $w^{\top} L=0$ and
$\sum_{i=1}^{m}w_{i}=1$. Moreover the symmetric part of $L^{\top}\diag(w)$ is positive semi-definite.
\end{lemma}

With $\mathbb{R}_-$, $\mathbb{R}_+$ and $\R_{\geqslant 0}$ we denote the sets of
negative, positive and nonnegative real numbers, respectively. The $i$th row and $j$th column of a matrix $M$ are denoted as $M_{i,\cdot}$ and $M_{\cdot,j}$, respectively. And for simplicity, let $M_{\cdot,j}^{\top}$ denote $(M_{\cdot,j})^\top$.
The vectors $e_1,e_2,\ldots,e_n$ denote the canonical basis of $\R^n$.

In the rest of this section we give some definitions and notations regarding
Filippov solutions.

Let $\funcRdRd$ be a map from $\R^n$ to $\R^n$,
and let $2^{\R^n}$ denotes the collection of all subsets of $\R^n$.
We define the \emph{Filippov set-valued map} of $\funcRdRd$, denoted
$\calF[\funcRdRd]:\R^n\rightarrow 2^{\R^n}$,  as
\begin{equation}\label{e:def_filippovset}
\calF[\funcRdRd](x):=
\bigcap_{\delta>0}\bigcap_{\mu(S)=0}\overline{\mathrm{co}}\{\funcRdRd(B(x,
\delta)\backslash S)\},
\end{equation}
where $B(x,\delta)$ is the open ball centered at $x$ with radius $\delta>0$,
$S$ is a subset of $\R^n$,
$\mu$ denotes the Lebesgue measure and $\overline{\mathrm{co}}$ denotes the convex closure.
If $X$ is continuous at $x$, then $ \calF[\funcRdRd](x)$ contains only the point $X(x)$.
There are some useful properties about the Filippov set-valued map.

\begin{lemma}[\cite{paden1987}]\label{p:calculus for Filippov}
Calculus for $\calF$.
	\begin{enumerate}[(i)]
		\item Assume that $f:\R^m\rightarrow\R^n$ is locally bounded. Then $\exists N_f\subset \R^m, \mu(N_f)=0$ such that $\forall N\subset\R^m, \mu(N)=0$,
		\begin{equation}
		\calF[f](x)=\mathrm{co}\{\lim_{i\rightarrow\infty} f(x_i)\mid x_i\rightarrow x, x_i\notin N_f\cup N \}.
		\end{equation}
		\item Assume that $f_j:\R^m\rightarrow \R^{n_j}, j=1,\ldots,N$ are locally bounded, then
		\begin{equation}
		\calF\big[ \bigtimes_{j=1}^N f_j \big](x) \subset \bigtimes_{j=1}^N\calF[f_j](x),
		\end{equation}
		where $\bigtimes$ represents the Cartesian product.
		\item Let $g:\R^m\rightarrow\R^{p\times n}$ be $C^0$ and $f:\R^m\rightarrow\R^n$ be locally bounded; then
		\begin{equation}
		\calF[gf](x)=g(x)\calF[f](x),
		\end{equation}
		where $gf(x):=g(x)f(x)\in\R^p$.
	\end{enumerate}
\end{lemma}

\begin{lemma}\label{p:fili_monotone}
	For an increasing function $\varphi:\R\rightarrow\R$, the Filippov set-valued map satisfies that
	\begin{enumerate}[(i)]
		\item  $\calF[\varphi](x)=[\varphi(x^-),\varphi(x^+)]$ where $\varphi(x^-),\varphi(x^+)$ are the left and right limit of $\varphi$ at $x$, respectively;
		\item  for any $x_1<x_2$, and $\nu_i\in\calF[\varphi](x_i), i=1,2,$ we have $\nu_1\leq\nu_2$.
	\end{enumerate}
\end{lemma}

\begin{pf}
	This can be seen as a straightforward deduction from Lemma \ref{p:calculus for Filippov} (i) and the definition of increasing functions.
\end{pf}

By using the fact that monotone functions are continuous almost everywhere, and the definition of right and left limits, we have following lemma.

\begin{lemma}\label{p:monotone}
	For an increasing function $\varphi:\R\rightarrow\R$,
	\begin{enumerate}[(i)]
		\item  $\calF[\varphi](x)=\{\varphi(x) \}$ for almost all $x$;
		\item  the right (left) limit, i.e., $\varphi(x^+)$ ($\varphi(x^-)$) is right (left) continuous for all $x$.
	\end{enumerate}
\end{lemma}

A \emph{Filippov solution} of the differential equation $\dot{x}(t)=\funcRdRd(x(t))$ on $[0,t_1]\subset\R$ is
an absolutely continuous function $x:[0,t_1]\rightarrow\R^n$ that
satisfies the differential inclusion
\begin{equation}\label{e:differential_inclusion}
\dot{x}(t)\in \calF[\funcRdRd](x(t)),
\end{equation}
for almost all $t\in[0,t_1]$.
A Filippov solution $t\mapsto x(t)$ is \emph{complete} if it is defined for all $t\in[0,\infty)$. Since the Filippov solutions of a discontinuous system \eqref{e:differential_inclusion} are not necessarily unique, we need to specify two types of invariant sets. A set $\calR\subset\R^n$ is called \emph{weakly invariant} for \eqref{e:differential_inclusion} if, for each $x_0\in \calR$, at least one complete solution of \eqref{e:differential_inclusion} with initial condition $x_0$ is contained in $\calR$. Similarly, $\calR\subset \R^n$ is called \emph{strongly invariant} for \eqref{e:differential_inclusion} if, for each $x_0\in \calR$, every complete solution of \eqref{e:differential_inclusion} with initial condition $x_0$ is contained in $\calR$. For more details, see \cite{cortes2008,filippov2013differential}.

Let $\funcRdR$ be a map from $\R^n$ to $\R$. The right directional derivative
of $\funcRdR$ at $x$ in the direction of $v\in \R^n$
is defined as
\begin{equation*}
\funcRdR'(x;v)=\lim_{h\rightarrow 0^+} \frac{\funcRdR(x+hv)-\funcRdR(x)}{h},
\end{equation*}
when this limit exists. The generalized derivative of $\funcRdR$ at $x$ in the
direction of $v\in \R^n$ is given by
\begin{equation*}
\begin{aligned}
\funcRdR^o(x;v) & =\limsup_{\begin{subarray}{c}
	y\rightarrow x  \\
	h\rightarrow 0^+
	\end{subarray}}
\frac{\funcRdR(y+hv)-\funcRdR(y)}{h} \\
& = \lim_{\begin{subarray}{c}
	\delta \rightarrow 0^+ \\
	\epsilon \rightarrow 0^+
	\end{subarray}}
\sup_{\begin{subarray}{c}
	y\in B(x,\delta)\\
	h\in[0,\epsilon)
	\end{subarray}}
\frac{\funcRdR(y+hv)-\funcRdR(y)}{h}.
\end{aligned}
\end{equation*}
We call the function $\funcRdR$ \emph{regular} at $x$ if $ \funcRdR'(x;v)$ and
$\funcRdR^o(x;v)$ are equal for all $v \in \R^n$. In particular, convex function is regular (see \cite{Clarke1990optimization}).

If $\funcRdR : \R^n \rightarrow \R$ is locally Lipschitz, then its {\it generalized gradient}
$\partial
\funcRdR:\R^n\rightarrow 2^{\R^n}$ is defined by
\begin{equation}
\partial \funcRdR(x):=\mathrm{co}\{\lim_{i\rightarrow\infty} \nabla
\funcRdR(x_i):x_i\rightarrow x, x_i\notin S\cup \Omega_{\funcRdR} \},
\end{equation}
where $\nabla$ denotes the gradient operator, $\Omega_{\funcRdR} \subset\R^n$ the set of points where
$\funcRdR$ fails to
be differentiable and $S\subset\R^n$ a set of Lebesgue measure zero that can be
arbitrarily
chosen to simplify the computation. The resulting set $\partial \funcRdR(x)$ is independent of the choice of $S$, see \cite{Clarke1990optimization}.

Given a  set-valued map $\calF:\R^n\rightarrow
2^{\R^n}$, the \emph{set-valued Lie derivative}
$\tilde{\mathcal{L}}_{\calF}\funcRdR:\R^n\rightarrow 2^{\R}$
of a locally Lipschitz function $\funcRdR:\R^n\rightarrow \R$  with respect to
$\calF$ at $x$ is
defined as
\begin{equation}\label{e:set-valuedLie}
\begin{aligned}
\tilde{\mathcal{L}}_{\calF}\funcRdR(x) := & \{a\in\R  \mid \textnormal{there
	exists } \nu\in\calF(x) \textnormal{ such that } \\
& \zeta^T\nu=a
\textnormal{
	for all } \zeta\in \partial \funcRdR(x)\}.
\end{aligned}
\end{equation}
If $\calF$ takes convex and compact values, then for each $x$, $\tilde{\mathcal{L}}_{\calF}\funcRdR(x)$ is a closed and bounded interval in $\R$, possibly empty.

The following result is a generalization of LaSalle's invariance principle to differential inclusions \eqref{e:differential_inclusion} with non-smooth
Lyapunov functions.
\begin{lemma}[LaSalle Invariance Principle, \cite{cortes2008}]\label{chap_preli:thm_stability}
Let $\funcRdR:\R^n\rightarrow\R$ be a locally Lipschitz and regular function. Let $S\subset \R^n$ be compact and strongly invariant for \eqref{e:differential_inclusion}, and assume that $\max \tilde{\calL}_{\calF[\funcRdRd]} f(y)\leq 0$ for each $y\in S$, where we define $\max\emptyset=-\infty$. Then, all solutions $x:[0,\infty)\rightarrow \R^n$ of \eqref{e:differential_inclusion} starting at $S$ converge to the largest weakly invariant set $M$ contained in
\begin{equation}
S\cap\overline{\{y\in\R^n\mid
	0\in\tilde{\calL}_{\mathcal{F}[\funcRdRd]}f(y)\}}.
\end{equation}
Moreover, if the set $M$ consists of a finite number of points, then the limit of each solution starting in $S$ exists and is an element of $M$.
\end{lemma}

At the end of this section, we list two potential Lyapunov functions.
\begin{lemma}[Prop. 2.2.6, Ex. 2.2.8 in \cite{Clarke1990optimization}]\label{regular_Lipschitz_lyapunov}
	The following functions are regular and Lipschitz continuous:
	\begin{equation}\label{e:VandW}
	V(x):= \max_{i\in\calI} x_i, \qquad W(x):= -\min_{i\in\calI} x_i.
	\end{equation}
\end{lemma}

\section{Multi-Agent Systems with Nonlinear Measurements}\label{s:measure}

In this section we consider a network of $n$ agents with a communication topology given by a weighted directed graph $\calG=(
\calV,\calE,A)$. Agent $i$ receives information from agent $j$ if and only if
there is an edge from node $v_j$ to node $v_i$ in the graph $\calG$. 
Consider the following nonlinear consensus protocol
\begin{equation}\label{e:nonlinear1}
\dot{x} = -Lf(x),
\end{equation}
where $f(x)=[f_1(x_1),\ldots,f_n(x_n)]^T$ and $f_i:\R\rightarrow\R$.
Throughout this paper, the following assumption is essential.

\begin{assumption}\label{as:monotone}
	The function $f_i:\R \rightarrow \R$ is an increasing function satisfying  $\lim_{x_i\rightarrow +\infty}f_i(x_i)>0$ and $\lim_{x_i\rightarrow -\infty}f_i(x_i)<0$.
\end{assumption}
Note that we do \emph{not} assume continuity of $f_i$. Examples of functions satisfying Assumption \ref{as:monotone} include sign function and quantizers. We understand the solution of \eqref{e:nonlinear1} in the Filippov sense, i.e., we consider the differential inclusion
\begin{equation}
	\begin{aligned}
	\dot{x} & \in \calF[-Lf(x)](x) \\
	&= -L\calF[f](x),
	\end{aligned}
\end{equation}
where the equality is implied by Lemma \ref{p:calculus for Filippov} (iii). Furthermore, by Lemma \ref{p:calculus for Filippov} (ii), the previous dynamical inclusion satisfies
\begin{equation}\label{e:nonlinear1_fili}
\begin{aligned}
\dot{x} \in -L \bigtimes_{i=1}^n\calF[f_i](x_i)
 := \calK_1(x).
\end{aligned}
\end{equation}
The existence of a Filippov solution can be guaranteed by the monotonicity of $f_i$, which indicates the local existence of solutions, see \cite{cortes2008}. Furthermore, we assume the complete solution of \eqref{e:nonlinear1_fili} exists for any initial condition.

Denote
\begin{equation}\label{e:D}\calD_1=\{ x\in\R^n \mid \exists a\in\R \textnormal{ s.t. } a\mathds{1}_n\in \bigtimes_{i=1}^n\calF[f_i](x_i)\}.
\end{equation}

\begin{lemma}
	Assumption \ref{as:monotone} holds, then set $\calD_1$ is closed.
\end{lemma}

\begin{pf}
	Take any sequence $\{y^k\}\subset\R ^n$ satisfying $\lim_{k\rightarrow\infty}y^k=x$ and $y^k\in\calD_1, k=1,2,\ldots$, we shall show that $x\in\calD_1$. Without loss of generality, we can assume the sequence $y^k_i$ converge to $x_i$ from one side, i.e., $y^k_i<x_i$ or $y^k_i>x_i$.
	
	Note that $y^k\in\calD_1$ implies that $\cap_{i=1}^n\calF[f_i](y^k_i)\neq \emptyset.$ For the case $y^k_i>x_i$, we have $f_i(y^{k-}_i)\geq f_i(x^-_i)$, $f_i(y^{k+}_i)\geq f_i(x^+_i)$ and $\lim_{k\rightarrow\infty}f_i(y^{k+}_i)= f_i(x^+_i)$ which is based on Lemma \ref{p:monotone} (ii). Hence we have
	\begin{equation}
	[\lim_{k\rightarrow\infty}f_i(y^{k-}_i),\lim_{k\rightarrow\infty}f_i(y^{k+}_i)]\subset [f_i(x^-_i),f_i(x^+_i)].
	\end{equation} Similarly, for the case $y^k_i<x_i$, this is also true. Then $\cap_{i=1}^n\calF[f_i](x_i)\neq \emptyset$, i.e., $x\in\calD_1$.
\end{pf}

\quad

\begin{theorem}\label{th:SC_nonlinear1}
	Suppose $\calG$ is a strongly connected digraph or connected undirected graph and Assumption \ref{as:monotone} holds. Then all Filippov solutions of \eqref{e:nonlinear1_fili} converge asymptotically to $\calD_1$.
\end{theorem}

\begin{pf}
	
	Consider the Lyapunov function $V_1(x)=w^TF(x)$ where $w\in\R^n_{+}$ is given by Lemma \ref{connected_laplacian} and \begin{equation*}
	F(x)=[F_1(x_1),\ldots,F_n(x_n)]
	\end{equation*} with $F_i(x_i)=\int_0^{x_i}f_i(\tau)d\tau$.
	It can be verified that $V_1\in\calC^0$ and $V_1$ is convex which implies that $V_1$ is regular. Moreover, by the monotonicity of $f_i$, we have $\partial F_i(x_i)=[f_i(x_i^-),f_i(x_i^+)]=\calF[f_i](x_i)$. Hence $V_1$ is locally Lipschitz continuous. Moreover, by Assumption \ref{as:monotone}, the function $V_1$ is radially unbounded. Indeed, $\lim_{x_i\rightarrow\infty}\int_{0}^{x_i}f_i(\tau)d\tau = \infty.$
	
	Let $\Psi_1$ be defined as
	\begin{equation}
	\Psi_1 = \{t\geq 0 \mid \textnormal{both } \dot{x}(t) \textnormal{ and } \frac{d}{dt}V_1(x(t)) \textnormal{ exist} \}.
	\end{equation}
	Since $x$ is absolutely continuous and $V_1$ is locally Lipschitz, we can let $\Psi_1=\R_{\geq 0}\setminus\bar{\Psi}_1$ where $\bar{\Psi}_1$ is a Lebesgue measure zero set. By Lemma 1 in \cite{Bacciotti1999}, we have
	\begin{equation}
	\frac{d}{dt}V_1(x(t))\in \tilde{\mathcal{L}}_{\calK_1}V_1(x(t)),
	\end{equation}
	for all $t\in\Psi_1$ and hence that the set $\tilde{\mathcal{L}}_{\calK_1}V_1(x(t))$ is nonempty for all $t\in\Psi_1$. For $t\in\bar{\Psi}_1$, we have that $\tilde{\mathcal{L}}_{\calK_1}V_1(x(t))$ is empty, and hence $\max \tilde{\mathcal{L}}_{\calK_1}V_1(x(t))< 0$. In the following, we only consider $t\in\Psi_1$. Moreover, in the proofs of the rest theorems in this paper, we always focus on a subset of $\R_{\geq 0}$ on which the set-valued Lie derivative of the corresponding Lyapunov functions are nonempty.
	
	The gradient of $V_1$ is given as
	\begin{equation}
	\partial V_1(x) =\co \{ \diag(w)\nu \mid \nu\in \bigtimes_{i=1}^n\calF[f_i](x_i) \}.
	\end{equation}
	Then $\forall a\in \tilde{\mathcal{L}}_{\calK_1}V_1(x(t))$, we have that $\exists u\in\bigtimes_{i=1}^n\calF[f_i](x_i)$ such that
	\begin{equation}
	a = -u^T L^T \diag(w)\nu
	\end{equation}
	for all $\nu \in\bigtimes_{i=1}^n\calF[f_i](x_i)$. A special case is that $\nu=u$, which implies that $a\leq 0$ by Lemma \ref{connected_laplacian}. Hence we have $\max \tilde{\mathcal{L}}_{\calK}V_1(x(t))\leq 0.$ Moreover, $a=0$ if and only if $\bigtimes_{i=1}^n\calF[f_i](x_i)\cap \spa\{\mathds{1}_n\}\neq\emptyset.$ Hence, by the fact that $\calD_1$ is closed, we have $\overline{\{x\in\R^n\mid 	0\in\tilde{\calL}_{\calK}V_1(x)\}}= \calD_1$.
	By Theorem \ref{chap_preli:thm_stability}, all the Filippov trajectories converges into the largest weakly invariant set containing in $\overline{\{x\in\R^n\mid 	0\in\tilde{\calL}_{\calK}V_1(x)\}}.$ Hence the conclusion holds.
\end{pf}

For homogenous systems, the requirement to graph $\calG$ can be weakened.

\begin{theorem}\label{th:ST_nonlinear1}
	Suppose $\calG$ is a digraph containing a spanning tree and the nonlinear functions in \eqref{e:nonlinear1} can be formulated as  $f(x)=[\bar{f}(x_1),\bar{f}(x_2),\ldots,\bar{f}(x_n)]$ where $\bar{f}$ satisfies Assumption \ref{as:monotone}. Then all Filippov solutions of \eqref{e:nonlinear1_fili} converge asymptotically to
	\begin{equation}\label{e:D_1}
	\calD_2=\{ x\in\R^n \mid \exists a\in\R \textnormal{ s.t. } a\mathds{1}_n\in \bigtimes_{i=1}^n\calF[\bar{f}](x_i)\}.
	\end{equation}
\end{theorem}

\begin{pf}
	In this case, the differential inclusion \eqref{e:nonlinear1_fili} can be written as
	\begin{equation}\label{e:nonlinear1_fili_barf}
	\begin{aligned}
	\dot{x}   \in -L \bigtimes_{i=1}^n\calF[\bar{f}](x_i)
	 := \calK_2(x).
	\end{aligned}
	\end{equation}
	
We divide the proof into five steps.
		
{\it (i)} Let's see the behaviors of the trajectories corresponding to roots. Noting the fact that the subgraph corresponding to the roots is strongly connected, by Theorem \ref{th:SC_nonlinear1}, all Filippov solutions of \eqref{e:nonlinear1_fili_barf} converge to
\begin{equation}
\{x\mid \exists a \textnormal{ s.t. } a\in\calF[\bar{f}](x_i), \forall i\in\calI_r \}.
\end{equation}
where $\calI_r =\{i\in\calI \mid v_i \text{ is a root of }\calG\}$.
	
{\it (ii)} Consider candidate Lyapunov functions $V$ as given in \eqref{e:VandW}.
Let $x(t)$ be a trajectory of \eqref{e:nonlinear1_fili_barf} and define
\begin{equation}\label{e:alpha}
\alpha(x(t))=\{k\in\calI\mid x_k(t)=V(x(t))\}.
\end{equation}
Denote $x_i(t)=\overline{x}(t)$ for $i\in\alpha(x(t))$.
The generalized gradient of $V$ is given as [\cite{Clarke1990optimization}, Example 2.2.8]
\begin{equation}
\partial V(x(t))  = \mathrm{co}\{e_k\in\R^n \mid k \in \alpha(x(t)) \}.
\end{equation}

Similar to the proof of Theorem \ref{th:SC_nonlinear1}, we can define $\Psi_2$ and we only consider $t\in\Psi_2$ such that $\tilde{\mathcal{L}}_{\calK_2}V(x(t))$ is nonempty and $\R_{\geq 0}\setminus\Psi_2$ is a Lebesgue measure zero set.
For $t \in \Psi_2$, let $a\in\tilde{\mathcal{L}}_{\calK_2}V(x(t))$.
By definition, there exists a $\nu^a \in \bigtimes_{i=1}^n\calF[\bar{f}](x_i)$ such that $a = (-L\nu^a)^{\top}\cdot \zeta$ for all
$\zeta\in\partial V(x(t))$. Consequently, by choosing $\zeta = e_k$ for $k \in \alpha(x(t))$, we observe that $\nu^a$ satisfies
	\begin{equation} \label{e:nu-alpha}
	   -L_{k,\cdot}	\nu^a  = a \qquad \forall k \in \alpha(x(t)).
	\end{equation}

Next, we want to show that $\max \tilde{\mathcal{L}}_{\calK_2}V(x(t))\leq 0$ for all $t\in\Psi_2$ by considering two possible cases: $\calI_r\nsubseteq\alpha(x(t))$ or
$\calI_r\subseteq\alpha(x(t))$.



If $\calI_r\subset \alpha(x(t))$, there are two subcases. First, $|\calI_r|=1$, i.e., there is only one root, denoted as $v_i$. Then $L_{i,\cdot}=0$, hence $L_{i,\cdot}\nu=0$ for any $\nu\in\bigtimes_{i=1}^n\calF[\bar{f}](x_i)$. By the observation \eqref{e:nu-alpha}, we have $\tilde{\mathcal{L}}_{\calK_2}V(x(t))=\{0 \}$. Second, $|\calI_r|\geq 2$. By the fact that the subgraph spanned by the roots is strongly connected, there exists $w_i>0$ for $i\in\calI_r$ such that $\sum_{i\in\calI_r}w_iL_{i,\cdot}=0_n, $
which implies that
\begin{equation}
\sum_{i\in\calI_r}w_iL_{i,\cdot}\nu=0
\end{equation}
for any $\nu\in\bigtimes_{i=1}^n\calF[\bar{f}](x_i)$. Again, by the observation \eqref{e:nu-alpha}, we have $\tilde{\mathcal{L}}_{\calK_2}V(x(t))=\{0 \}$.

If $\calI_r \nsubseteq \alpha(x(t))$, i.e., there exists $i\in\calI_r\setminus \alpha(x(t))$. We define a subset $\alpha'(\nu)$ as
\begin{equation}
\alpha'(\nu)=\{i\in\alpha(x(t)) \mid \nu_i=\max_{i\in\alpha(x(t))} \nu_i \}
\end{equation}
for any $\nu\in\bigtimes_{i=1}^n\calF[\bar{f}](x_i)$. From Lemma \ref{p:fili_monotone} (ii),  for any $j\in\alpha'(\nu)$, we know that $\nu_j=\max \nu_i$, thus $L_{j,\cdot}\nu\geq 0$. By the fact that the choice of $\nu$ is arbitrary in $\bigtimes_{i=1}^n\calF[\bar{f}](x_i)$ and the observation \eqref{e:nu-alpha}, we have $\tilde{\mathcal{L}}_{\calK_2}V(x(t))\subset \R_{\leq 0}$. Moreover, denoting
\begin{equation}
\calE_{\alpha(x)}=\{e_{ij}\in\calE\mid j\in\alpha(x) \},
\end{equation}
we shall show that $0\in\tilde{\mathcal{L}}_{\calK_2}V(x)$ if and only if $\exists \nu\in \bigtimes_{i=1}^n\calF[\bar{f}](x_i)$ such that $\nu_i=\nu_j$ for any $e_{ij}\in\calE_{\alpha(x)}$, which is equivalent to  $\calF[\bar{f}](x_i)\cap\calF[\bar{f}](x_j)\neq \emptyset$ for all $e_{ij}\in\calE_{\alpha(x)}$.
The sufficient part is straightforward, in fact we can take $\nu_i=\nu_j=f(\overline{x}^-)$ for any $e_{ij}\in\calE_{\alpha(x)}$. Then $0\in\tilde{\mathcal{L}}_{\calK_2}V(x)$. The necessary part can be proved as follows. Since $0\in\tilde{\mathcal{L}}_{\calK_2}V(x)$, there exists $\nu\in \bigtimes_{i=1}^n\calF[\bar{f}](x_i)$ such that $L_{j,\cdot}\nu= 0$ for any $j\in\alpha(x)$. Then this $\nu$ satisfies that $\alpha'(\nu)=\alpha(x)$. Indeed, if $\alpha'(\nu) \subsetneqq \alpha(x)$, then for any $j\in\alpha'(\nu)$ with $e_{ij}\in\calE$ and $i\notin\alpha'(\nu)$, $L_{j,\cdot}\nu<0$. 
Hence $\alpha'(\nu)=\alpha(x)$. Furthermore, by using the same argument, we have for any $e_{ij}\in\calE$ satisfying $i\notin\alpha(x)$ and $j\in\alpha(x)$, $f(\overline{x}^-)\in\calF[\bar{f}](x_i)$.


{\it (iii)} For the Lyapunov functions $W$ as given in \eqref{e:VandW},
denote
\begin{equation}\label{e:beta}
\beta(x(t))=\{i\in\calI\mid x_i(t)=-W(x(t))\},
\end{equation}
and
$x_i(t)=\underline{x}(t)$ for $i\in\beta(x(t))$, and $\calE_{\beta(x(t))} =\{e_{ij}\in\calE\mid j\in\beta(x(t)) \}$. By using similar computations, we find that $\max \tilde{\mathcal{L}}_{\calK_2}W(x(t))\leq 0$ and $0\in\tilde{\mathcal{L}}_{\calK_2}W(x(t))$ if and only if $\exists \nu\in \bigtimes_{i=1}^n\calF[\bar{f}](x_i)$ such that $\nu_i=\nu_j$ for any $e_{ij}\in\calE_{\beta(x(t))}$, which is equivalent to  $\calF[\bar{f}](x_i)\cap\calF[\bar{f}](x_j)\neq \emptyset$ for all $e_{ij}\in\calE_{\beta(x(t))}$.

{\it (iv)} So far we have that $V(x(t))$ and $W(x(t))$ are not increasing along the trajectories $x(t)$ of the
system \eqref{e:nonlinear1_fili_barf}. Hence, the trajectories are bounded and remain in the set $[\underline{x}(0),\overline{x}(0)]^n$ for all $t\geq0$.
Therefore, for any $N\in\R_+$, the set $S_N=\{x\in\R^n \mid \|x\|_{\infty}\leqslant N\}$ is
strongly invariant for \eqref{e:nonlinear1_fili_barf}.
By Theorem \ref{chap_preli:thm_stability}, we
have that all solutions of \eqref{e:nonlinear1_fili_barf} starting in $S_N$
converge to
the largest weakly invariant set $M$ contained in
\begin{equation}\label{e: 0inLVandLW}
\begin{aligned}
S_N & \cap\overline{\{x\in\R^n:
	0\in\tilde{\mathcal{L}}_{\calK_2}V(x)\}} \\
	& \cap\overline{\{x\in\R^n:
	0\in\tilde{\mathcal{L}}_{\calK_2}W(x)\}}.
\end{aligned}
\end{equation}

{\it (v)} We have proved the asymptotic stability of the system. Next we will prove that the set $\calD_2$ is strongly invariant and for any $x_0\notin\calD_2$, all the solution satisfying $x(0)=x_0$ will converge to $\calD_2$.

We start with the strong invariance of $\calD_2$. Notice that by the monotonicity of $\bar{f}$ we can reformulate $\calD_2$ as
\begin{equation}
\calD_2=\{x\mid \calF[\bar{f}](\underline{x})\cap \calF[\bar{f}](\overline{x})\neq \emptyset \}.
\end{equation}
For any $x_0\in\calD_2$, we have known that any trajectories starting from $x_0$, $V(x(t))$ and $W(x(t))$ are not increasing. Hence $\overline{x}(t)\leq \overline{x}_0$ and $\underline{x}(t)\geq \underline{x}_0$ for all $t\geq 0$ which, by Lemma \ref{p:fili_monotone}, implies that $\calF[\bar{f}](\underline{x}(t))\cap \calF[\bar{f}](\overline{x}(t))\neq \emptyset$ for all $t$ and $x(t)$ satisfying $x(0)=x_0$. Then $x(t)\in\calD_2$ which implies that $\calD_2$ is strongly invariant.

Next we show that for any $x_0\notin\calD_2$, all the solution satisfying $x(0)=x_0$ will converge to $\calD_2$. We will prove it by contradictions. Indeed, we assume that there exists $x_0\notin\calD_2$ and one solution $\tilde{x}(t)$ satisfying $\tilde{x}(0)=x_0$ does not converge to $\calD_2$. Since the set $\calD_2$ is strongly invariant, we have $\tilde{x}(t)\notin\calD_2$ for all $t\geq 0.$ Then $\calF[\bar{f}](\underline{\tilde{x}})\cap \calF[\bar{f}](\overline{\tilde{x}})= \emptyset$, where
\begin{equation*}
\begin{aligned}
\overline{\tilde{x}}  = \lim_{t\rightarrow\infty}V(\tilde{x}(t)),~
\underline{\tilde{x}}  = -\lim_{t\rightarrow\infty}W(\tilde{x}(t)).
\end{aligned}
\end{equation*}
Hence there exists a constant $C>0$, such that
\begin{equation}
d(\calF[\bar{f}](\underline{\tilde{x}}), \calF[\bar{f}](\overline{\tilde{x}}))>C
\end{equation}
where $d(S_1,S_2)=\inf_{y_1\in S_1,y_2\in S_2}d(y_1,y_2)$ is the \emph{distance} between two sets $S_1$ and $S_2$. For any $i,j\in\calI$ with $i\neq j$, there exists a vector $w^{ij}\in\R^n$ such that $w^{ij^\top} L=(e_i-e_j)^T$. For each pair $i,j\in\calI$, we choose one $w^{ij}$ and collect all the $w^{ij}$ for $i,j\in\calI$ in the set $\Omega$. Notice that there are only finite number of vectors in $\Omega$.  Then for any $t, i\in\alpha(\tilde{x}(t))$ and $j\in\beta(\tilde{x}(t))$, we have $\overline{\tilde{x}}(t)\geq\overline{\tilde{x}}$ and $ \underline{\tilde{x}}(t)\leq\underline{\tilde{x}} $. Moreover, since $\tilde{x}(t)$ is uniformly bounded, there exist a constant $\tau$ which does not depend on $t$ such that for any $s\in[t, t+\tau]$
\begin{equation}\label{e:differential_lower_bound}
w(s)^T \dot{x}(s)>\frac{C}{2}.
\end{equation}
where $w:\R\rightarrow\Omega$ is piecewise constant and $w(s)=w^{ij}$ with $i\in\alpha(t),j\in\beta(t)$ for $s\in[t,t+\tau]$. Note that for any $T$, the function $w(s)^T\dot{x}(s)$ is Lebesgue integrable on $[0, T]$, and by \eqref{e:differential_lower_bound} we have
\begin{equation}
\int_{0}^T w(s)^\top \dot{x}(s)ds>\frac{C}{2}T
\end{equation}
which converge to infinity as $T\rightarrow \infty$.
This is a contradiction to the fact that $w(s)$ is globally bounded and for any $T<\infty$ and $i\in\calI$, $\int_{0}^T \dot{x}_i(s)ds$ is bounded. Hence we have for any $x_0\notin\calD_2$, all the solution satisfying $x(0)=x_0$ will converge to $\calD_2$. Here ends the proof.

\end{pf}

\begin{remark}\label{r:positive_nonlinear1}
	From the proof of Theorem \ref{th:ST_nonlinear1}, we know the maximal components of the trajectories of the system \eqref{e:nonlinear1_fili_barf} are not increasing while the minimal ones are not decreasing. Hence \eqref{e:nonlinear1_fili_barf} is a positive system (see e.g., \cite{Rantzer2011}), i.e., with positive initial conditions, the trajectories will be positive for all the time. However, the system \eqref{e:nonlinear1_fili} is in general not a positive system.
\end{remark}

\begin{remark}
	The stability of system \eqref{e:nonlinear1} under more general assumptions than the ones in Theorem \ref{th:ST_nonlinear1}, namely the nonlinear functions $f_i$ are different for each agent and the underlying graph is directed which contains a spanning tree,  is still an open problem.
\end{remark}

\section{Applications to quantized consensus }\label{s:quantize}

In this section, we shall apply the results in the previous section to the quantized multi-agent systems. There are three types of  quantizers, namely the symmetric, asymmetric and logarithmic quantizer mainly considered in the literature
	\begin{align}
	q^s(z) & = \Big\lfloor \frac{z}{\Delta}+\frac{1}{2} \Big\rfloor \Delta,\nonumber\\
	q^a(z) & = \Big\lfloor \frac{z}{\Delta} \Big\rfloor \Delta,\\
	q^l(z) & = \begin{cases}
	\sign(z)\exp\Big(\mathtt{q}^s\big(\ln(|z|)\big)\Big) & \textrm{ if } z\neq 0,\\
	0 & \textrm{ if } z=0,
	\end{cases} \nonumber
	\end{align}
respectively.

There are some properties about these quantizers. First, for the symmetric quantizer $\mathtt{q}^s$ we have: {\it (i)}
$|\mathtt{q}^s(z)-z|\leq \frac{\Delta}{2}$; {\it (ii)} $\mathtt{q}^s(z)=-\mathtt{q}^s(-z)$. Second, for the asymmetric quantizer $\mathtt{q}^a$, the following relation holds: $0\leq z-\mathtt{q}^a(z)\leq \Delta$. Finally, the logarithmic quantizer $\mathtt{q}^l$ satisfies: {\it (i)} $\mathtt{q}^l(z)=-\mathtt{q}^l(-z)$; {\it (ii)} $|\mathtt{q}^l(z)-z|<\big(\exp(\frac{\Delta}{2})-1 \big)|z|$.

By denoting $q(x)=(q_1(x_1),\ldots,q_n(x_n)^T$ where $q_i:\R\rightarrow\R, i=1,\ldots,n$ is a quantizer, the system \eqref{e:nonlinear1} can be written as
\begin{equation}\label{e:quantizedconsensus_general_1}
\dot{x} = -Lq(x).
\end{equation}

For the case of digraphs, we consider the quantizers satisfy that $q_i=q^s, \forall i\in\calI$ and the system \eqref{e:quantizedconsensus_general_1} can be written as
\begin{equation}\label{e:quantizedconsensus1}
\dot{x} = -Lq^s(x).
\end{equation}
In this case the set $\calD_2$ defined as \eqref{e:D_1} is given as
\begin{equation}
\{x\in\R^n \mid \exists k\in\Z \textnormal{ such that } k\Delta\mathds{1}_n\in\mathcal{F}[q^s](x)\},
\end{equation}
which is equivalent to
\begin{align}
\calQ:=& \{x\in\R^n \mid \exists k\in\Z \textnormal{ s. t. } \\& (k-\frac{1}{2})\Delta \leq x_i\nonumber \leq (k+\frac{1}{2})\Delta, \forall i \in\calI \}.
\end{align}
It is known that without the precise measurement of the states, exact consensus can not be achieved in principle. Instead, the notation of \emph{practical consensus} will be employed. We say that the state variables of the agents converge to \emph{practical consensus}, if $x(t)\rightarrow\calQ$ as $t\rightarrow\infty$.

Based on Theorem \ref{th:ST_nonlinear1}, we have the following results which is an extension of the result in Section 3 of \cite{Ceragioli2011}. More precisely, we generalize the result in \cite{Ceragioli2011} to the digraphs containing a spanning tree.

\begin{corollary}\label{c:quantized_measure}
     Suppose $\calG$ is a digraph containing a spanning tree. Then all Filippov solutions of \eqref{e:quantizedconsensus1} converge asymptotically to
    practical consensus, i.e., $\calQ$.
\end{corollary}

\begin{remark}
	By Proposition 1 in \cite{Ceragioliphdthesis}, the Krasovskii and Filippov solutions of \eqref{e:quantizedconsensus1} are equivalent. Hence Corollary \ref{c:quantized_measure} holds for all Krasovskii solutions as well.
\end{remark}

\begin{remark}
When the underlying topology is a strongly connected digraph or connected undirected graph, Theorem \ref{th:SC_nonlinear1} implies stability of the hybrid quantized system where  agents can have different quantizers, i.e.,
\begin{equation}
\dot{x} = -Lq^*(x),
\end{equation}
where $q^*_i$ can be $q^s, q^a$ or $q^l$.
\end{remark}

\section{Conclusions}\label{s:conclusion}

In this paper, we considered a general nonlinear consensus protocol, namely the multi-agent systems with nonlinear measurement of their states. Here we assumed the nonlinear functions to be monotonic increasing without any continuity constraints. The solutions of the dynamical systems were understood in the sense of Filippov. We proved asymptotic stability of the systems defined on different topologies. More precisely, we considered the systems defined on undirected graphs or digraphs containing a spanning tree. Finally, we applied the results to quantized consensus. Future interesting problems include the switching topology and robustness to uncertainties.

\bibliography{references}             








\appendix

\end{document}